%% file: DPLQGMain.tex
\documentclass[letterpaper, 10 pt, conference, draftcls, onecolumn]{ieeeconf}  

\IEEEoverridecommandlockouts                              

\usepackage{graphicx} 
\usepackage{amsmath} 
\usepackage{amssymb}  
\usepackage{color}
\usepackage[linesnumbered, ruled, vlined]{algorithm2e}
\usepackage{mathtools}
\mathtoolsset{showonlyrefs}
\graphicspath{{.}}
\newtheorem{dfn}{Definition}
\newtheorem{lemma}{Lemma}

\newtheorem{problem}{Problem}
\newtheorem{theorem}{Theorem}
\newtheorem{remark}{Remark}
\newtheorem{assumption}{Assumption}

\newcommand{\R}{\mathbb{R}}
\newcommand{\N}{\mathbb{N}}
\newcommand{\E}{\mathbb{E}}
\newcommand{\logdet}{\log\det}
\newcommand{\barsig}{\overline{\Sigma}}
\newcommand{\tell}{\tilde{\ell}}
\newcommand{\diag}{\textnormal{diag}}
\newcommand{\adj}{\textnormal{Adj}}
\newcommand{\sys}{\mathcal{G}}
\newcommand{\mech}{\mathcal{M}}
\renewcommand{\Pr}{\mathbb{P}}

\newcommand{\sval}{s}
\newcommand{\dpvar}{\sigma}
\newcommand{\tfin}{K_f}

\title{\LARGE \bf
Privacy in Feedback: The Differentially Private LQG
}

\author{Matthew Hale, Austin Jones, Kevin Leahy
\thanks{Matthew Hale is with the Department of  Mechanical and Aerospace Engineering at the University of Florida, Gainesville, FL USA. 
Austin Jones and Kevin Leahy are with the Massachusetts Institute of Technology Lincoln Laboratory, Lexington, MA USA}%
\thanks{DISTRIBUTION STATEMENT A. Approved for public release: distribution unlimited.}%
\thanks{This material is based upon work supported under Air Force Contract No. FA8721-05-C-0002 and/or FA8702-15-D-0001. Any opinions, findings, conclusions or recommendations expressed in this material are those of the author(s) and do not necessarily reflect the views of the U.S. Air Force.}%
}%

\begin{document}

\maketitle

	\input{abstract}
\input{intro}

\input{PrivacyBackground}

\input{Problem}

\input{Results}
\input{CaseStudy}

\input{Conclusions}


\bibliographystyle{plain}{}
\bibliography{sources}

\end{document}

%% file: abstract.tex
\begin{abstract}
Information communicated within cyber-physical systems (CPSs) is often used
in determining the physical states of such systems, and malicious
adversaries may intercept these communications in order to infer
future states of a CPS or its components. Accordingly, there arises
a need to protect the state values 
of a system. Recently, the notion of differential privacy has been used to protect state
trajectories in dynamical systems, and it is this notion of privacy
that we use here to protect the state trajectories of CPSs. 
We incorporate a cloud computer to coordinate the agents comprising
the CPSs of interest, and the cloud offers the ability to remotely coordinate
many agents, rapidly perform computations, and broadcast the results, making
it a natural fit for systems with many interacting agents or components. 
Striving for broad applicability, we solve infinite-horizon 
linear-quadratic-regulator (LQR) problems, and each agent protects its own state trajectory
by adding noise to its states before they are sent to the cloud. 
The cloud then uses these state values to generate optimal inputs for the
agents. As a result,
private data is fed into feedback loops at each iteration, and each noisy term
affects every future state of every agent. In this paper, we show that
the differentially private LQR problem can be related to the well-studied
linear-quadratic-Gaussian (LQG) problem, and we provide bounds on
how agents' privacy requirements affect the cloud's ability to generate
optimal feedback control values for the agents. These results are illustrated
in numerical simulations. 

\end{abstract}

%% file: intro.tex
\section{Introduction}

Many distributed cyber-physical systems (CPSs), such as robotic swarms and smart power networks, involve  generating and transmitting data that could be exploited by a malicious adversary, such
as the position of a robot or power usage of an individual. 
Ideally, this data would be secure, meaning that an adversary would not be able to access any meaningful information.  
In practice, secure data transmission and consumption may not be practical due to physical limitations, e.g., limited hardware power, 
or because some information must be shared to attain a system's primary objective, e.g., ensuring that a power grid's production
is not less than its consumption. In such cases, information must be shared in a manner which makes it useful to the
intended recipient while protecting it from any eavesdroppers. 

Recently, privacy of this form has been achieved using differential privacy. Differential privacy sanitizes data 
with differentially private mechanisms (DPMs), which are randomized functions of data whose distributions over outputs are
close for similar input data. Practically, this means that changes in a single agent's input data are not likely to be
seen at a DPM's output, protecting individuals' sensitive data while still allowing for accurate conclusions to be drawn
about groups of agents. 
For example, DPMs may be used in  a micro-grid to sanitize pricing information in such a way that participants may approximately 
optimize the price they pay for electricity, while a malicious agent with access to the same data cannot recover a single individual's power usage profile \cite{hale2015cloud,doe2010}.   
Since DPMs are randomized functions, their usage may degrade the performance of a system in comparison to a system without privacy simply because
noise is added where it would otherwise be absent. Despite this drawback, DPMs are suitable for protecting sensitive data because they are: 
\begin{enumerate}
\item Robust to arbitrary side-information, meaning that additional information about the data-producing entities does not yield much additional 
information about the true value of privatized data \cite{kasiviswanathan2008thesemantics}
\item Robust to post-processing, meaning that an adversary cannot increase their knowledge of the true value of the privatized data by 
applying a post-hoc transformation to a DPM's output \cite{dwork2014algorithmic} 
\item Robust to mechanism knowledge, meaning that an adversary with complete knowledge of the mechanism used to privatize data has no 
advantage over an adversary without this knowledge \cite{dwork2014algorithmic,dwork2006calibrating}. 
\end{enumerate}

It is natural to ask what the utility of DPMs is given the existence of encryption techniques.  DPMs may be used along with encryption in a layered protection strategy 
in which encryption provides intrusion prevention and differential privacy provides intrusion mitigation.  Further, encryption may not be desirable for data transmission 
in  mobile ad-hoc networks due to the additional computation, power consumption, and latency that encryption introduces \cite{zhang2014lightweight,sahingoz2014networking,yang2004security},
whereas using a DPM requires only minor changes to implement. 

Over the past decade, researchers have developed a rich theory of differential privacy for database protection that has translated into practical applications 
\cite{dwork2014algorithmic,greenberg2016apple,machanavajjhala2008privacy}.  In addition, recent results have extended these works to protect entire trajectories produced by entities, 
such as dynamical systems or networks of agents, that evolve over time in a predictable manner \cite{leny14}. 
More information on the relevant aspects of differential privacy may be found in Section~\ref{sec:privacy_review}.  
The trajectory privatization strategies promulgated in some existing works have assumed that the end use of the 
privatized trajectories is to perform \em post-hoc \em statistical analyses \cite{leny14b,leny14}.   
However, many CPSs that publish or share potentially sensitive data, such as traffic networks or teams of autonomous systems, 
make decisions based on past data that will affect the state of the system in the future.  
This suggests the possibility that the performance degradation induced by applying DPMs to data transmitted across the 
network may compound over time to continuously degrade performance. While privacy has been
used in other feedback configurations, e.g., in \cite{wang14} for database-style privacy with
completely correlated noise over time, 
to the best of our knowledge, these feedback effects have not been previously addressed for trajectory-level privacy across infinite time horizons. 

In this paper, as a first step towards characterizing the effects of trajectory-level privacy in feedback loops, we analyze a 
multi-agent system in which each agent transmits its state to a central cloud computer \cite{hale14} 
that 
then computes optimal inputs for these agents pointwise in time. These inputs are transmitted back to the agents which apply 
them in their local state updates, and then this process of exchanging information and updating states repeats. 
We consider a discrete time linear quadratic regulation 
problem in which each system evolves according to linear dynamics and the cloud computes inputs to minimize 
the expected average quadratic state and input costs incurred at each iteration of the system's evolution. We further quantify the effects of privacy by bounding the entropy seen by the cloud
in this setup. 
If the true value 
of the state information were reported over the network, then an adversary could use that information to trivially recover the state 
trajectories of the individual agents. Instead, our implementation provides strong theoretical privacy guarantees to each
agent while still allowing for optimal control values to be computed at each timestep. 
The contribution of this paper therefore
consists of the differentially private LQG control algorithm along with the quantification of privacy's impact upon the system. 

The rest of the paper is organized as follows. Section~\ref{sec:privacy_review} reviews the relevant aspects of differential
privacy and provides the privacy background needed for the rest of the paper. 
Then, Section~\ref{sec:problem} defines the cloud-enabled differentially private LQR problem in which state measurements are perturbed
by noise. 
Section \ref{sec:results} then gives a closed-form solution to this problem and provides an algorithm for solving it on the
cloud-based architecture we use. Then we present an upper bound for how the noise introduced by privacy affects 
the certainty in the cloud's estimate of the joint state of the network. Next, we provide some interpretation of these results in 
numerical simulations in Section \ref{sec:case}, and then Section~\ref{sec:conclusions} provides conclusions and directions for future
work.

%% file: PrivacyBackground.tex
\section{Review of Differential Privacy} \label{sec:privacy_review}
In this section we review  necessary background 
on differential privacy and its implementation for
dynamical systems. A survey on differential privacy for databases can be
found in \cite{dwork2014algorithmic}, and the adaptation of differential
privacy to dynamical systems is developed in \cite{leny14}, which
forms the basis for the exposition in this section.  
We first review the basic definitions required for
differential privacy and then describe the privacy mechanism
we use in the remainder of the paper.  We use the notation $[\ell] = \{1, \ldots, \ell\}$ for $\ell \in \N$.

\subsection{Differentially Private dynamical Systems}
Fundamentally, differential privacy is used to
protect the data of individual users while still allowing
for accurate statistical analyses of groups of users \cite{dwork2006calibrating}.
In this work, a user's data is his or her state
trajectory, and we therefore keep these trajectories
differentially private, thereby protecting individual users' activities,
while still allowing for meaningful computations to be performed
on collections of these trajectories. 

To develop our privacy implementation,
suppose there are $N$ agents, with agent $i$'s state trajectory denoted
$x_i$, and the $k^{th}$ element of this trajectory
denoted $x_i(k) \in \R^{n_i}$ for some $n_i \in \N$. 
The state trajectory
$x_i$ is contained in the set $\tell^{n_i}_2$, namely
the set of sequences of vectors in $\R^{n_i}$ 
whose finite truncations all have finite $\ell_2$-norm.
Formally, we define the truncation operator $P_T$
over trajectories according to
\begin{equation}
P_T[x] = \begin{cases} x(k) & k \leq T \\
                      0   & k > T 
       \end{cases},
\end{equation}
and we say that $x_i \in \tell^{n_i}_2$ if and only if $P_T[x_i] \in \ell^{n_i}_2$
for all $T \in \N$. The trajectory space for all $N$ agents in aggregate
is then $\tell^{n}_2 = \tell^{n_1}_2 \times \cdots \times \tell^{n_N}_2$,
where $n = \sum_{i \in [N]} n_i$. 

Differential privacy for dynamical systems seeks to make
\emph{adjacent} state trajectories produce outputs which are similar
in a precise sense to be defined below. An adversary or eavesdropper
only has access to these privatized outputs and does not have access to the
underlying state trajectories which produced them. 
Therefore, the definition of adjacency
is a key component of any privacy implementation as it
specifies which pieces of sensitive data must be made approximately
indistinguishable to any adversary or eavesdropper in order to
protect their exact values. 
We define our adjacency relation now over the space $\tell^n_2$ defined
above. 

\begin{dfn} \label{def:adj} (\emph{Adjacency}) Fix $B > 0$. 
The adjacency relation $\adj_B$ is defined
for all $v, w \in \tell^{n}_2$ as
\begin{equation}
\adj_B(v, w) = \begin{cases}
        1 & \parbox[t]{.3\textwidth}{$\|v - w\|_{\ell_2} \leq B$} \\
        0 & \text{otherwise}. \tag*{$\triangle$}
        \end{cases}
\end{equation}
\end{dfn}
In words, two state trajectories are adjacent if and only if 
the $\ell_2$ distance between them is less than or equal to $B$. 
This adjacency relation therefore requires that every agent's state trajectory
be made approximately indistinguishable from all other state trajectories
not more than distance $B$ away. 

The choice of adjacency relation gives rise to the notion of sensitivity
of a dynamical system, which we define now.

\begin{dfn} \emph{(Sensitivity)} 
The $p$-norm sensitivity of a dynamical system $\sys$ is 
the greatest distance between two output trajectories
which correspond to adjacent state trajectories, i.e., 
\begin{equation}
\Delta_p\sys := \sup_{x, \tilde{x} \mid \adj_B(x, \tilde{x})=1} \|\sys(x) - \sys(\tilde{x})\|_{\ell_2}. \tag*{$\triangle$}
\end{equation}
\end{dfn}

Below, the sensitivity of a system will be used in determining how much noise
must be added to make it differentially private. 
Noise is added by a \emph{mechanism}, which is a randomized map
used to enforce differential privacy. 
For private dynamical systems,
it is the role of a mechanism
to take sensitive trajectories as inputs and approximate their corresponding
outputs with stochastic processes. 
Below, we give a formal definition of
differential privacy for dynamical systems which specifies
the probabilistic guarantees that a mechanism must provide. 
To do so, we first fix a probability space $(\Omega, \mathcal{F}, \mathbb{P})$.
This definition
considers outputs in the space $\tell^q_2$ and uses 
a {$\sigma$-algebra} over $\tell^q_2$, denoted $\Sigma^q_2$.
Details on the formal construction of this $\sigma$-algebra 
can be found in \cite{leny14}.

\begin{dfn} \label{def:dp} (\emph{Differential Privacy}, \cite{leny14}) Let 
$\epsilon > 0$ and $\delta \in (0, 1/2)$ be given. 
A mechanism ${M : \tell^n_2 \times \Omega \to \tell^q_2}$ is $(\epsilon, \delta)$-differentially private
if and only if, for all adjacent $x, \tilde{x} \in \tell^n_2$, we have
\begin{equation}
\Pr[M(x) \in S] \leq e^{\epsilon}\Pr[M(\tilde{x}) \in S] + \delta \textnormal{ for all } S \in \Sigma^q_2. \tag*{$\triangle$}
\end{equation}
\end{dfn}

A key feature of this definition is that one may tune the extent to which
$x$ and $\tilde{x}$ are approximately indistinguishable. The values
of $\epsilon$ and $\delta$ are user-specified and
together determine the level of privacy provided
by the mechanism $M$, with smaller values of each generally providing
stronger privacy guarantees to each user. These parameters
are commonly selected so that $\epsilon \in [0.1, \log 3]$ and $\delta \in [0, 0.05]$  \cite{leny14}.
In the literature, $(\epsilon, 0)$-differential privacy is called 
 $\epsilon$-differential privacy, while referring to
$(\epsilon, \delta)$-differential privacy is understood to imply
$\delta > 0$. In this paper, we consider only cases
in which $\delta > 0$ because enforcing $(\epsilon, \delta)$-differential privacy
can be done via the Gaussian mechanism.

\subsection{The Gaussian Mechanism}
The Gaussian mechanism is defined in terms of the $\mathcal{Q}$-function, which is given by
\begin{equation}
\mathcal{Q}(y) = \frac{1}{\sqrt{2\pi}}\int_{y}^{\infty} e^{-\frac{z^2}{2}}dz. 
\end{equation}

We now formally define the Gaussian mechanism, which adds Gaussian noise
point-wise in time to the output of a dynamical system to make it private. 

\begin{dfn} \label{def:gaussmech} \emph{(Gaussian Mechanism; \cite{leny14}, Theorem 3)}
Let privacy parameters $\epsilon > 0$ and $\delta \in (0, 1/2)$ be given. 
Let $\sys$ denote a dynamical system with state trajectories in $\tell^n_2$ and 
outputs in $\tell^q_2$, and denote its $2$-norm sensitivity by $\Delta_2 \sys$. 
Then the Gaussian mechanism for $(\epsilon, \delta)$-differential privacy
takes the form $\mech(x) = \sys(x) + w$,
where $w(k) \sim \mathcal{N}(0, \sigma I_q)$, $I_q$ is the $q \times q$
identity matrix, and $\sigma$ is bounded according to
\begin{equation}
\sigma \geq \frac{\Delta_2\sys}{2\epsilon}\big(K_{\delta} + \sqrt{K_{\delta}^2 + 2\epsilon}\big), \textnormal{ with } K_{\delta} := \mathcal{Q}^{-1}(\delta). \tag*{$\triangle$}
\end{equation}
\end{dfn}

For convenience, we define the function
\begin{equation}
\kappa(\delta, \epsilon) = \frac{1}{2\epsilon}\big(K_{\delta} + \sqrt{K_{\delta}^2 + 2\epsilon}\big),
\end{equation}
and we will use the Gaussian mechanism to enforce differential privacy for
the remainder of the paper. 

%% file: Problem.tex
\section{Problem Formulation}
\label{sec:problem}
In this section we define the multi-agent LQ problem to be solved privately.
First, we give the agents' dynamics and equivalent network-level dynamics.
Then, we define the cloud-based control scheme and give a problem
statement that requires optimizing running state and control costs while protecting
agents' privacy. 

We use the notation $P \succeq 0$ and $P \succ 0$ to indicate that a matrix $P$ is positive 
semi-definite or positive definite matrices, respectively. 
Similarly, we write $A \succeq B$ and $A \succ B$ to indicate
that $A - B \succeq 0$ and $A - B \succ 0$, respectively. 
We use the notation
\begin{equation}
\diag(P_1, \ldots, P_n) := P_1 \oplus \cdots \oplus P_n
\end{equation}
to denote the matrix direct sum of $P_1$ through $P_n$. 
If these matrices are square, 
the notation $\diag(P_1, \ldots, P_n)$ denotes 
a block diagonal matrix
with blocks $P_1$ through $P_n$; if the elements $P_i$ are scalars,
this is simply a diagonal matrix with these scalars on its main diagonal. 

\subsection{Multi-Agent LQ Formulation}
As in Section~\ref{sec:privacy_review}, consider a collection of $N$ agents, indexed over
$i \in [N]$. 
At time $k$, agent $i$ has state 
vector $x_i(k) \in \R^{n_i}$ for some $n_i \in \N$, which evolves according to
\begin{equation}
x_i(k+1) = A_ix_i(k) + B_iu_i(k) + w_i(k),
\end{equation}
where $u_i(k) \in \R^{m_i}$, $m_i \in \N$, is agent $i$'s input at time $k$, 
$w_i(k) \in \R^{n_i}$ is process noise, and where
$A_i \in \R^{n_i \times n_i}$ and $B_i \in \R^{n_i \times m_i}$ are
time-invariant matrices. The process noise is distributed as $w_i(k)\sim\mathcal{N}(0, W_i)$,such that the $w_i$ are uncorrelated and such that $0 \prec W_i \in \R^{n_i \times n_i}$. 
 
To analyze the agents' behavior in aggregate, 
we define the network-level state and control vectors as
\begin{equation}
x(k) = \left(\begin{array}{c} x_1(k) \\ \vdots \\ x_N(k)\end{array}\right) \in \R^{n} \textnormal{ and } u(k) = \left(\begin{array}{c} u_1(k) \\ \vdots \\ u_N(k) \end{array}\right) \in \R^{m},
\end{equation}
where $n = \sum_{i \in [N]} n_i$ and $m = \sum_{i \in [N]} m_i$. 
By further defining
\begin{equation}
w(k) = \left(\begin{array}{c} w_1(k) \\ \vdots \\ w_n(k)\end{array}\right),
\end{equation}
$A = \diag(A_1, \ldots, A_N) \in \R^{n \times n}$, and $B = \diag(B_1, \ldots, B_N) \in \R^{m \times m}$, we have the network-level state dynamics
\begin{equation} \label{eq:dynamics}
x(k+1) = Ax(k) + Bu(k) + w(k).
\end{equation}

We consider infinite-horizon problems in which the agents are jointly incur the  average cost
\begin{equation} \label{eq:cost}
J(x, u) = \lim_{\tfin \to \infty} \frac{1}{\tfin} \left[\sum_{k=1}^{\tfin} x(k)^TQx(k) + u(k)^TRu(k)\right],
\end{equation}
where $Q \in \R^{n \times n}$ and $R \in \R^{m \times m}$. 
The linear dynamics in Equation~\eqref{eq:dynamics} and quadratic cost in Equation~\eqref{eq:cost}
together define a linear-quadratic (LQ) optimal control problem, a canonical
class of optimal control problems \cite{anderson90,bertsekas05,sontag13}.
The problem considered below is distinct from the existing literature because it incorporates the 
exchange of differentially privatized data in order to protect agents' state trajectories, which, to the best of our
knowledge, has not been done before in an LQ problem.  More formally, we address Problem \ref{prob:main}.

\begin{problem} \label{prob:main}
Given an initial measurement $x(0) \in \R^{n}$ of the network state, 
\begin{equation}
\min_{u} J(x, u) = \lim_{\tfin \to \infty} \frac{1}{\tfin} \left[\sum_{k=1}^{\tfin} x(k)^TQx(k) + u(k)^TRu(k)\right]
\end{equation}
subject to the network-level dynamics
\begin{equation}
x(k+1) = Ax(k) + Bu(k) + w(k)
\end{equation}
while keeping each agent's state trajectory differentially private with agent-provided
 privacy parameters $(\epsilon_i,\delta_i)$. \hfill $\blacklozenge$
\end{problem}

\subsection{Cloud-Based Aggregation} \label{ss:cloud_description}
In Problem~\ref{prob:main}, the cost $J$ is generally non-separable due to its
quadratic form This non-separability means that $J$ cannot be minimized
by each agent using only knowledge of its own state. 
Instead, at each point in time, each
agent's optimal control vector will, in general, depend upon all states
in the network. To compute the required optimal control vectors, one
may be tempted to have an agent aggregate all agents' states and then
perform the required computations, though
this is undesirable for two reasons. First, doing so would inequitably
distribute computations in the network by making a single agent
gather all information and perform all computations on it. In an 
application with battery-powered agents, this results
in the designated aggregating agent's lifetime being significantly shorter than the others.
Second, in some cases, the cost $J$ may be known
only by a central aggregator and not by the agents themselves,
rendering this approach unusable. 

Due to the need for aggregation of network-level information and computations
performed upon it, we introduce a cloud computer into the system.
The technology of cloud computing provides the ability to communicate
with many agents, rapidly perform demanding computations, and broadcast the
results. Moreover, these capabilities can be seamlessly and remotely added
into networks of agents, making the cloud a natural choice of central
aggregator for this work \cite{hale14}. 

 At time $k$, the cloud requests
from agent $i$ the value $y_i(k) = C_ix_i(k)$, where $C_i \in \R^{n_i \times n_i}$
is a constant matrix. 
We refer to $y_i$ as either agent $i$'s output or as agent $i$'s
state measurement. 
To protect its state trajectory,
agent $i$ sends a differentially private form of $y_i(k)$ to the cloud.
Using these values from each agent, the cloud computes 
the optimal controller for each agent at time $k$. The cloud sends
these control values to the agents, the agents use them in their local state updates,
and then this process of exchanging information repeats. 

Next we specify the exact manner in which
each agent makes its transmissions to the cloud differentially private.
Following that, Section~\ref{sec:results} provides a complete solution algorithm
for Problem~\ref{prob:main}.

\subsection{Differentially Private Communications}
At time $k$, agent $i$'s transmission of $y_i(k)$ to the cloud can potentially
reveal its state trajectory, thereby compromising agent $i$'s privacy. 
Accordingly, agent $i$ adds noise to the state measurement $y_i(k)$ to keep
its state trajectory differentially private. Agent $i$ first selects its own privacy
parameters ${\epsilon_i > 0}$ and $\delta_i \in (0, 1/2)$. With agent $i$
using the adjacency relation $\adj_{b_i}$ for some $b_i > 0$ (cf. Definition~\ref{def:adj}),
we have the following elementary lemma concerning the $2$-norm sensitivity
of $y_i$, which we denote $\Delta_2 y_i$.

\begin{lemma} \label{lem:ysensitivity}
The $2$-norm sensitivity of $y_i$ satisfies
${\Delta_2 y_i \leq \sval_{1}(C_i)b_i}$,
where $\sval_{1}(\cdot)$ denotes the maximum singular value of a matrix. 
\end{lemma}
\emph{Proof:}
Consider two trajectories $x_i, \tilde{x}_i \in \tell^{n_i}_2$ such that
$\adj_{b_i}(x_i, \tilde{x}_i) = 1$, and denote $y_i(k) = C_ix_i(k)$ and
$\tilde{y}_i(k) = C_i\tilde{x}_i(k)$. Then we find
\begin{equation}
\|y_i - \tilde{y}_i\|_{\ell_2} \leq \|C_i\|_2\|x_i - \tilde{x}_i\|_{\ell^2} \leq \sval_{1}(C_i)b_i. \tag*{$\blacksquare$}
\end{equation}

%
%
With the bound in Lemma~\ref{lem:ysensitivity}, agent $i$ sends the cloud the private output
\begin{equation}
\bar{y}_i(k) := y_i(k) + v_i(k) = C_ix_i(k) + v_i(k),
\end{equation}
where $v_i(k) \sim \mathcal{N}(0, \dpvar_iI_{n_i})$ is a Gaussian
random variable, and where $\dpvar_i \geq \kappa(\delta_i, \epsilon_i)\sval_{1}(C_i)b_i$
in accordance with Definition~\ref{def:gaussmech}.  
Defining the matrix $V = \diag(\dpvar_1I_{n_1}, \ldots, \dpvar_NI_{n_N})$, the aggregate
privacy vector in the system at time $k$ is
\begin{equation}
v(k) = (v_1(k)^T \,\, \cdots \,\, v_N(k)^T)^T
\end{equation}
Defining the matrix
\begin{equation}
V = \diag(\dpvar_1I_{n_1}, \ldots, \dpvar_NI_{n_N}),
\end{equation}
we have $v(k) \sim \mathcal{N}(0, V)$ for all $k \in \N$. 

In this privacy implementation, we assume that all matrices $A_i$, $B_i$, and $C_i$ are public information.
In addition, we impose the following assumption on the dynamics and cost that are considered.

\begin{assumption} \label{as:matrices}
In the cost $J$, $Q \succ 0$ and $R \succ 0$. In addition, the pair $(A, B)$ is controllable,
and there exists a matrix $F$ such that $C = F^TF$ and such that the pair $(A, F)$ is observable.
\hfill $\triangle$
\end{assumption}

\begin{figure}
\begin{center}
\includegraphics[width=0.7\columnwidth]{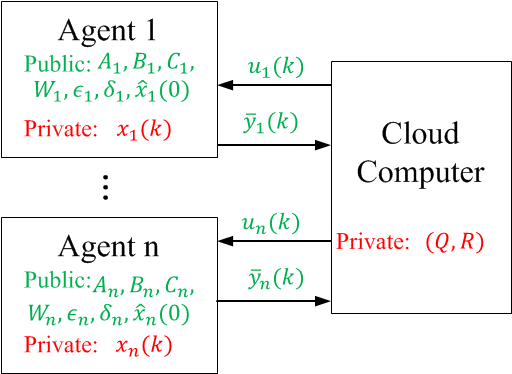} 
\caption{\label{fig:block} Summary of information available over network.}
\end{center}
\end{figure}

In addition, we assume that $\hat{x}_i(0)$, the expected value of agent $i$'s
initial state, is publicly known, along with agent $i$'s privacy parameters
$\epsilon_i$ and $\delta_i$. 
We also assume that the matrices $Q$ and $R$ are known only to the cloud
and that the cloud does not share these matrices with any agents or eavesdroppers.  The information
available to agents, the cloud, and eavesdroppers is summarized in Figure \ref{fig:block} 
With this privacy implementation in hand, we next
give a reformulation of Problem~\ref{prob:main}.

\subsection{Reformulation of Problem~\ref{prob:main}} \label{ss:newprob}
Having specified exactly the sense in which privacy is implemented in Problem~\ref{prob:main},
we now give an equivalent statement of Problem~\ref{prob:main} as an LQG control problem \cite{athans71}. 

\begin{problem} \label{prob:newprob}
Minimize 
\begin{equation} 
J(x, u) = \lim_{\tfin \to \infty} \frac{1}{\tfin} \left[\sum_{k=1}^{\tfin} x(k)^TQx(k) + u(k)^TRu(k)\right]
\end{equation}
over all $u$ with $u(k) \in \R^{m}$, subject to the dynamics and output equations
\begin{align}
x(k+1) &= Ax(k) + Bu(k) + w(k) \\
y(k) &= Cx(k) + v(k),
\end{align}
where $C = \diag(C_1, \ldots, C_N)$ and where $\hat{x}(0) := \E[x(0)]$ is known and each agent has specified its privacy parameters $(\epsilon_i,\delta_i)$. 
\hfill $\blacklozenge$
\end{problem}

%% file: Results.tex
\section{Results: Private LQG Control} \label{sec:results}
In this section we solve Problem~\ref{prob:newprob}, with the understanding
that this solution also solves Problem~\ref{prob:main}. Problem~\ref{prob:newprob}
takes the form of a Linear Quadratic Gaussian (LQG) problem, which is well-studied 
in the control literature \cite{athans71,astrom70}, We first solve Problem~\ref{prob:newprob}
using established techniques and then quantify the effects of privacy
upon the system. Below, we use the notation $\lambda_i(\cdot)$ to denote the $i^{th}$ largest
eigenvalue of a matrix and the notation $\sval_i(\cdot)$ to denote the $i^{th}$ largest
singular value of a matrix. 

\subsection{Solving Problem~\ref{prob:newprob}}
Due to the process noise and privacy noise in Problem~\ref{prob:newprob}, the controllers we develop
cannot rely on the exact value of $x(k)$. Instead, the controllers used
rely on its expected
value given all past inputs and outputs, denoted
\begin{equation}
\hat{x}(k) := \E[x(k) \mid u(0), \ldots, u(k-1), y(0), \ldots, y(k-1)].
\end{equation}
Problem~\ref{prob:newprob} is an infinite-horizon discrete-time LQG problem,
and it is known \cite[Section 5.2]{bertsekas05} that the optimal controller
for such problems relies on $\hat{x}(k)$ in the form
\begin{equation} \label{eq:ustar}
u^*(k) = L\hat{x}(k).
\end{equation}
Here, 
\begin{equation} \label{eq:findL}
L = -(R + B^TKB)^{-1}B^TKA,
\end{equation}
and $K$ is the unique positive semidefinite solution to the
discrete algebraic Riccati equation
\begin{equation} \label{eq:findK}
K = A^TKA - A^TKB(R + B^TKB)^{-1}B^TKA + Q.
\end{equation}

Computing the state estimate $\hat{x}(k)$ can be done for infinite time horizons
using the time-invariant Kalman filter \cite[Section 5.2]{bertsekas05}
\begin{equation} \label{eq:kalmanfilter}
\hat{x}(k+1) + \barsig C^TV^{-1}\big(\bar{y}(k+1) - C(A + BL)\hat{x}(k)\big),
\end{equation}
where $\barsig$ is given by
\begin{equation} \label{eq:findbarsig}
\barsig = \Sigma - \Sigma C^T(C\Sigma C^T + V)^{-1}C\Sigma,
\end{equation}
and where $\Sigma$ is the unique positive semidefinite solution to the discrete
algebraic Riccati equation
\begin{equation} \label{eq:findsig}
\Sigma = A\Sigma A^T - A\Sigma C^T(C\Sigma C^T + V)^{-1}C\Sigma A^T + W.
\end{equation}

As discussed in Section~\ref{ss:cloud_description}, the cloud is responsible for
generating control values for each agent at each point in time. Therefore, the cloud
will run the Kalman filter in Equation~\eqref{eq:kalmanfilter} and will likewise generate
control values using Equation~\eqref{eq:ustar}. Thus the flow of information 
in the network first has agent $i$ send $\bar{y}_i(k)$ to the cloud at time $k$,
and then has the cloud form $\bar{y}(k) = \big(\bar{y}_1(k)^T, \ldots, \bar{y}_N(k)^T\big)^T$ and use
it to estimate $\hat{x}(k+1)$. Next, the cloud computes $u^*(k) = L\hat{x}(k)$,
and finally the cloud sends $u^*_i(k)$ to agent $i$. We emphasize that the
cloud only sends $u^*_i(k)$ to agent $i$ and, because sharing this control
value can threaten its privacy, agent $i$ does not share this control value
with any other agent. 

With respect to implementation, the matrices $K$, $L$, $\barsig$, and $\Sigma$ can all
be computed offline by the cloud a single time before the network begins solving
Problem~\ref{prob:main}. 
Then these matrices can be used repeatedly in generating
control values to send to the agents, substantially reducing the cloud's
computational burden at runtime and allowing it to quickly generate state
estimates and control values. With this in mind, 
we now state the full solution to Problem~\ref{prob:main} below
in Algorithm~\ref{alg:main}. 

\begin{algorithm} \label{alg:main}
 \KwData{Public information $A_i$, $B_i$, $C_i$, $\epsilon_i$, $\delta_i$, $\hat{x}_i(0)$}
 Initialize the cloud with all public information \\
 For all $i$, initialize agent $i$ with $A_i$, $B_i$, $C_i$. Agent $i$ chooses $(\epsilon_i, \delta_i)$ and
 knows its exact initial state $x_i(0)$ but does not share it. \\
 In the cloud, precompute $L$ using Equation~\eqref{eq:findL} and $\barsig$ according to Equation~\eqref{eq:findbarsig} \\
 \For{$k=0,1,2,\ldots$}{
   \For{$i=1,\ldots,N$}{
     Agent $i$ generates $v_i(k) \sim \mathcal{N}(0, \dpvar_iI_{n_i})$ \\
     Agent $i$ sends the cloud the private output $\bar{y}_i(k) := C_ix_i(k) + v_i(k)$
   }
   In the cloud, compute $\hat{x}(k)$ with Equation~\eqref{eq:kalmanfilter} \\
   In the cloud, compute $u^*(k) = L\hat{x}(k)$ and send $u^*_i(k)$ to agent $i$ \\
   \For{$i=1,\ldots,N$}{
     Agent $i$ updates its state via $x_i(k+1) = A_ix_i(k) + B_iu^*_i(k) + w_i(k)$
   }
 }
 \caption{Differentially Private LQG (Solution to Problem~\ref{prob:main})}
\end{algorithm}

\subsection{Quantifying the Effects of Privacy}
Algorithm~\ref{alg:main} solves Problem~\ref{prob:newprob}, though it is intuitively clear
that adding noise for the sake of privacy will diminish the ability of the cloud
to compute optimal control values. Indeed, the purpose of differential privacy in this problem
is to protect an agent's state value from the cloud (and any eavesdroppers). Here, we quantify the trade-off between privacy and the ability
of the cloud to accurately estimate the agents' state values.

The cloud both implements a Kalman filter and computes the controller $u^*(k)$, though
noise added for privacy only affects the Kalman filter. Note that
the network-level optimal controller defined in Equation~\eqref{eq:ustar} takes
the form $u^*(k) = L\hat{x}(k)$, 
while the optimal controller with no privacy noise (i.e., $v(k) \equiv 0$ for all $k$)
would take the same form with the expectation carried out only over the process
noise. Similarly, in the deterministic case of $v(k) \equiv 0$ and $w(k) \equiv 0$ for all
$k$, the optimal controller is
$u^*(k) = Lx(k)$.
The optimal controller in all three
cases is linear and defined by $L$, and the only difference between these
controllers comes from changes in the estimate of $x(k)$ or in using $x(k)$ exactly.
This is an example of the so-called ``certainty equivalence principle''
\cite{bertsekas05}, and it demonstrates that the controller in Algorithm~\ref{alg:main}
is entirely unaware of the presence or absence of noise. Instead, the effects of noise
in the system are compensated for by the Kalman filter defined in Equation~\eqref{eq:kalmanfilter},
and we therefore examine the impact of privacy noise upon the Kalman filter that informs
the controller used in Algorithm~\ref{alg:main}. Related work in \cite{leny14} specifies procedures
for designing differentially private Kalman filters, though here we use a Kalman filter to process
private data and quantify the effects of privacy upon the filter itself. 
To the best of our knowledge this is the first study to do so.

One natural way to study this trade-off is through bounding the information content in the Kalman
filter as a function of the privacy noise, as differential privacy explicitly seeks to mask sensitive information by corrupting
it with noise.  In our case, we consider the differential entropy of $\hat{x}$ as a proxy for the information
content of the signal \cite{cover2012elements}.  Shannon entropy has been used to quantify the effects of differential privacy in other
settings, including in distributed linear control systems where database-type
differential privacy is applied \cite{wang14}  Though differential entropy does not have the same axiomatic 
foundation as Shannon entropy, it is useful for Gaussian distributions because it bounds the extent of the sublevel sets of $\mathcal{R}^{-1}$ where  $\mathcal{R}(y) = 1-2\mathcal{Q}(y)$, i.e., the volume of covariance ellipsoids.  We can therefore quantify the effects of differentially private masking upon
the cloud by studying how privacy noise changes the $\logdet({\Sigma})$, which is within a constant factor of the differential entropy of $\hat{x}$. The analysis we present
differs from previous work because we apply trajectory-level differential privacy
and quantify the resulting differential entropy in the Kalman filter. 
Toward presenting this result, we have the following lemma concerning the determinants of matrices
related by a matrix inequality. 

\begin{lemma} \label{lem:logdet}
Let $P \succeq 0$ and $Q \succeq 0$ satisfy
$P \succeq Q$. Then $\textnormal{det}(P) \geq \textnormal{det}(Q)$ and
$\logdet(P) \geq \logdet(Q)$. 
\end{lemma}
\emph{Proof:} 
This follows from the monotonicity of $\log$ and the fact that $\lambda_i(X) \geq \lambda_i(Y)$ \cite[Theorem 16.F.1]{marshall09} whenever
$X \succeq Y$. 
\hfill $\blacksquare$

Next we have the following lemma that relates the determinant of a matrix to its trace.

\begin{lemma} \label{lem:amgm}
Let~$M \in \R^{n \times n}$ be a positive semi-definite matrix. Then
\begin{equation*}
\det(M) \leq \left(\frac{\textnormal{tr}(M)}{n}\right)^n.
\end{equation*}
\end{lemma}
\emph{Proof:} The positive semi-definiteness of~$M$ implies that all of its
eigenvalues are non-negative. The result then follows by applying
the AM-GM inequality to the eigenvalues of~$M$. \hfill $\blacksquare$

Finally, we have the following lemma which gives a matrix upper bound on solutions to
a discrete algebraic Riccati equation. 

\begin{lemma} \label{lem:sigbound}
Suppose $\Sigma$ is the unique positive semi-definite solution to the equation
\begin{equation} \label{eq:dare}
\Sigma = A(\Sigma - \Sigma C^T(C \Sigma C^T + V)^{-1}C \Sigma)A^T + W,
\end{equation}
and define
\begin{equation}
\gamma_i = \frac{\sigma_i^2W_{ii}}{\sigma_i^2 + C_{ii}^2W_{ii}},
\end{equation}
with $\Gamma = \diag(\gamma_1, \ldots, \gamma_n)$.
Then
\begin{equation}
\lambda_1(\Sigma) \geq s_n^2(A)\lambda_1(\Gamma) + \lambda_n(W) =: \eta.
\end{equation}
Furthermore, if 
\begin{equation}
s_1^2(A) < 1 + \big(s_n^2(A)\lambda_1(\Gamma) + \lambda_n(W)\big)\cdot \min_{i \in [n]} \frac{C_{ii}^2}{\sigma_i^2},
\end{equation}
then 
\begin{equation}
\Sigma \preceq \frac{\lambda_1(W)}{1 + \eta \lambda_n(C^TV^{-1}C) - s_1^2(A)}AA^T + W. 
\end{equation}
\end{lemma}
\emph{Proof:} This follows from 
Theorem~2 in \cite{lee97}. 
\hfill $\blacksquare$

We now present our main results on bounding the log-determinant of the steady-state 
Kalman filter error covariance in the private LQG problem. In particular, we study the 
log-determinant of the \emph{a priori} error covariance, $\Sigma$, because 
$\Sigma$ measures the error in the cloud's predictions about the state of the
system. Bounding the log-determinant of $\Sigma$ therefore bounds the
entropy seen at the cloud when making predictions about the state of the
network in the presence of noise added for privacy.

\begin{theorem} \label{thm:entropy}
Let $\Sigma$ denote the steady-state \emph{a priori} error covariance in the cloud's
Kalman filter in Algorithm~\ref{alg:main}. If all hypotheses of Lemma~\ref{lem:sigbound} hold, 
$C_i$ is square and diagonal for all $i$, and Assumption~\ref{as:matrices} holds, then
\begin{equation}
\logdet \Sigma < \left(\frac{\lambda_1(W)}{1 + \eta \lambda_n(C^TV^{-1}C) - s_1^2(A)}\right) \sum_{i=1}^{n} s_i^2(A) + \textnormal{tr}(W),
\end{equation}
with $\eta$ defined as in Lemma~\ref{lem:sigbound}, and
where the term
\begin{equation}
\eta\lambda_n(C^TV^{-1}C) = \sval_n^2(A)\cdot \max_{i \in [n]}\frac{\sigma_i^2 W_{ii}}{\sigma_i^2 + C_{ii}^2W_{ii}}\cdot \min_{i \in [n]} \frac{C_{ii}^2}{\sigma_i^2}
+ \lambda_n(W)\cdot \min_{i \in [n]} \frac{C_{ii}^2}{\sigma_i^2} 
\end{equation}
is the only one affected by privacy. 
\end{theorem}
\emph{Proof:} With the assumption that each $C_i$ is square and diagonal, we see that
\begin{equation}
C^TV^{-1}C = \textnormal{diag}\left(\frac{C_{11}^2}{\sigma_1^2}, \ldots, \frac{C_{nn}^2}{\sigma_n^2}\right),
\end{equation}
which gives
\begin{equation}
\lambda_{n}\Big(C^TV^{-1}C\Big) = \min_{1 \leq i \leq n} \frac{C_{ii}^2}{\sigma_i^2},
\end{equation}
and the form of~$\eta\lambda_n(C^TV^{-1}C)$ follows. 

Lemma~\ref{lem:sigbound} then gives
\begin{equation} \label{eq:main_ub}
\Sigma \preceq \frac{\lambda_1(W)}{1 + \eta\lambda_n(C^TV^{-1}C) - s_1^2(A)}AA^T + W,
\end{equation}
where applying Lemma~\ref{lem:logdet} gives
\begin{equation}
\det \Sigma \leq \det\bigg[\frac{\lambda_1(W)}{1 + \eta\lambda_n(C^TV^{-1}C) - s_1^2(A)}AA^T + W\bigg].
\end{equation}
The argument of the determinant on the right-hand side is a positive-definite matrix, and applying
Lemma~\ref{lem:amgm} therefore gives
\begin{equation}
\det \Sigma \leq \left(\frac{\textnormal{tr}\left(\frac{\lambda_1(W)}{1 + \eta\lambda_n(C^TV^{-1}C) - s_1^2(A)}AA^T + W\right)}{n}\right)^n. 
\end{equation}
Taking the natural logarithm of both sides, we find that
\begin{align}
\log\det\Sigma &\leq \log \left(\frac{\textnormal{tr}\left(\frac{\lambda_1(W)}{1 + \eta\lambda_n(C^TV^{-1}C) - s_1^2(A)}AA^T + W\right)}{n}\right)^n \\
             &= n \log \left(\frac{\textnormal{tr}\left(\frac{\lambda_1(W)}{1 + \eta\lambda_n(C^TV^{-1}C) - s_1^2(A)}AA^T + W\right)}{n}\right) \\
             &< n  \left(\frac{\textnormal{tr}\left(\frac{\lambda_1(W)}{1 + \eta\lambda_n(C^TV^{-1}C) - s_1^2(A)}AA^T + W\right)}{n}\right),
\end{align}
which follows from the fact that~$\log x < x$ for all~$x > 0$. Continuing, we find that
\begin{align}
\log \det \Sigma &< \textnormal{tr}\left(\frac{\lambda_1(W)}{1 + \eta\lambda_n(C^TV^{-1}C) - s_1^2(A)}AA^T + W\right) \\
                    &= \frac{\lambda_1(W)}{1 + \eta\lambda_n(C^TV^{-1}C) - s_1^2(A)} \textnormal{tr}(AA^T) + \textnormal{tr}(W),
\end{align}
which follows from the linearity of the trace operator. 
The theorem follows by expanding the trace of~$AA^T$ in terms of the singular values of~$A$. 
\hfill $\blacksquare$

\begin{remark}
To interpret Theorem~\ref{thm:entropy}, we consider the case in which $C_i = I_{n_i}$ for
all $i$, indicating that the cloud simply requests each agent's state at each point in time,
in which $W_i = \omega I_{n_i}$, indicating that each agent has identical noise covariance
for its process noise, and in which $\epsilon_i = \epsilon$ and $\delta_i = \delta$ for 
all $i$ and some $\epsilon > 0$ and $\delta \in (0, 1/2)$. In this case, $\sigma_i = \sigma$
for all $i$ and we find
\begin{equation}
\eta\lambda_n(C^TV^{-1}C) = \sval_n^2(A)\frac{\omega}{\sigma^2 + \omega} + \frac{\omega}{\sigma^2}. 
\end{equation}
If the other terms in the denominator can be neglected in Theorem~\ref{thm:entropy}, then one
finds that 
\begin{equation}
\logdet \Sigma \lessapprox \left(\frac{s_n^2(A)}{\sigma^2 + \omega} + \frac{1}{\sigma^2}\right)^{-1}\sum_{i=1}^{n}s_i^2(A) + n\omega,
\end{equation}
indicating that, in this case, a linear change in the variance of privacy noise
$\sigma$ incurs a roughly quadratic increase in entropy. \hfill $\lozenge$
\end{remark}

%% file: CaseStudy.tex
\section{Case Study}
\label{sec:case}
In this section, we use a two-agent example to illustrate the use of Algorithm~\ref{alg:main}. Consider a two-agent system 
with $x_1(k), x_2(k) \in \R^2$ for all $k$, and 
where each agent has identical dynamics and process noise covariance matrices

\begin{equation} \label{eq:dynMat}
A_i = \begin{bmatrix} 1 & 0.1 \\0 &  1 \end{bmatrix}, B_i = \begin{bmatrix}0 \\ 1 \end{bmatrix}, W_i = \begin{bmatrix} 1 & 0.5 \\ 0.5 & 1.0 \end{bmatrix}.
\end{equation}

In this case, we can interpret the first coordinate of each agent's position and the second coordinate as its velocity.  Consider the case in which the two agents select different privacy parameters, namely, 
\begin{equation} 
(\epsilon_1,\delta_1)= (0.1,0.01) \textnormal{ and } (\epsilon_2,\delta_2) = (1,0.5), 
\end{equation}
giving $\sigma_1 = 23.48$ and $\sigma_2 = 0.71$. 

\begin{figure}
\begin{center}
\includegraphics[width=0.75\columnwidth]{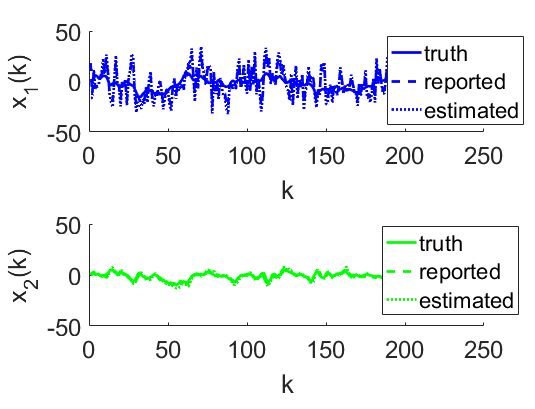} \\ (a) \\
\includegraphics[width=0.75\columnwidth]{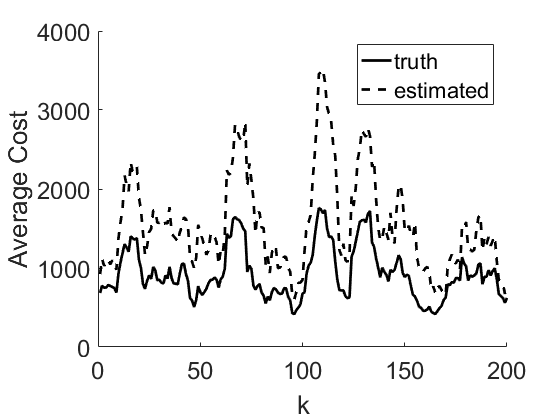} \\ (b)
\caption{\label{fig:traces} Examples of applying the DPLQG to a two agent network. (a) Agent  position trajectories (b) Moving average
cost.}
\end{center}
\end{figure}

These parameters mean that Agent 1 has chosen to be much more private than Agent 2. Algorithm~\ref{alg:main} was simulated for
this problem for $200$ timesteps, and the results of this simulation are shown in 
Figure~\ref{fig:traces}. For these simulations, the $Q,R$ matrices were randomly chosen positive definite matrices for which the 
off-diagonal elements were non-zero, meaning that the cost $J$ is not separable over agents. Note that the estimated trajectory of Agent 1, 
which is the optimal estimate (in the least squares sense) for both the aggregator and an adversary, shows substantial temporal variation 
and is in general very different from the true trajectory. However, even with the degraded information reported to the aggregator, the state 
of Agent 1 is still controlled closely to its desired value. And, although the aggregator substantially overestimates the incurred cost due 
to the injected privacy noise, the true incurred average cost over time remains bounded here. 

\begin{figure}
\begin{center}
\includegraphics[width=0.7\columnwidth]{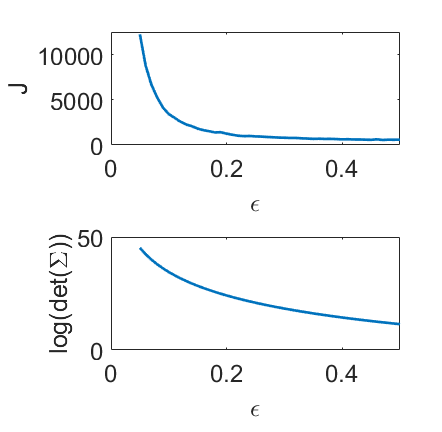} 
\caption{\label{fig:eps} Effect of privacy parameter $\epsilon$ on cost and estimation accuracy.}
\end{center}
\end{figure}

The effect of privacy parameters is further illustrated in Figure~\ref{fig:eps}. In this case, we consider a network with 4 agents with dynamics matrices as stated in Equation~\eqref{eq:dynMat},
and fix the parameters $\epsilon_i = \epsilon$ and $\delta_i = 0.25$ for all agents. 
Figure \ref{fig:eps} shows the effect on the cost $J$ (estimated from 2500 time steps) and the state estimate uncertainty $\logdet(\Sigma)$ 
as we vary the parameter $\epsilon$. These plots show that each quantity varies monotonically and smoothly with $\epsilon$, demonstrating that changing the overall
privacy level in the system indeed results in predictable changes in the behavior of the system, as predicted by Theorem~\ref{thm:entropy}.

%% file: Conclusions.tex
\section{Conclusions} \label{sec:conclusions}
In this paper, we have studied the problem of 
distributed linear quadratic optimal control with
agent-specified differential privacy requirements. To our knowledge,
this is the first paper to consider the effects of trajectory-level
differential privacy in feedback. Here, we have related this
problem to the well-studied LQG problem and have shown
how to bound the uncertainty in the network-level joint state
estimate in terms of the privacy parameters specified by the individual agents.
Future work includes
considering more general classes of control and estimation models and nested 
and sequential feedback loops.